\documentclass[hidelinks,11pt,a4paper,notitlepage,oneside]{amsart}
\usepackage[utf8]{inputenc}
\usepackage[english]{babel}

\usepackage{csquotes}
\usepackage[style=alphabetic,sorting=nyt,firstinits=true,maxnames=5,maxalphanames=5,backend=biber]{biblatex}
\usepackage{hyperref}

\usepackage{amsmath}
\usepackage{amsthm}
\usepackage{amssymb}
\usepackage{mathtools}
\usepackage{stmaryrd}
\usepackage{aligned-overset}
\usepackage{enumitem}
\usepackage{thmtools}
\usepackage{cleveref}

\theoremstyle{plain}
\newtheorem{theorem}{Theorem}[section]
\newtheorem{lemma}[theorem]{Lemma}
\newtheorem{proposition}[theorem]{Proposition}
\newtheorem{corollary}[theorem]{Corollary}

\theoremstyle{definition}

\newtheorem{remark}[theorem]{Remark}

\newcommand{\NN}{\mathbb{N}}
\newcommand{\RR}{\mathbb{R}}
\newcommand{\ZZ}{\mathbb{Z}}
\newcommand{\QQ}{\mathbb{Q}}

\newcommand{\FF}{\mathbb{F}}

\newcommand{\Prm}{\mathcal{P}}
\newcommand{\M}{\mathrm{M}}

\newcommand{\Sym}{\mathrm{Sym}}

\newcommand{\Dm}{\mathcal{D}}

\DeclareMathOperator{\id}{id}

\DeclareMathOperator{\Aut}{Aut}

\DeclareMathOperator{\spn}{span}

\DeclareMathOperator{\im}{im}

\DeclareMathOperator{\rk}{rk}

\newcommand{\OG}{\mathrm{O}}
\newcommand{\OP}{\mathrm{O}^+}
\newcommand{\OU}{\mathrm{O}^\prime}
\newcommand{\SOG}{\mathrm{SO}}

\newcommand{\SOU}{\mathrm{SO}^\prime}
\newcommand{\OT}{\widetilde{\mathrm{O}}}

\DeclareMathOperator{\E}{E}

\newcommand{\ZP}{{Z^{\mathrm{prim}}}}

\begin{document}

\title[]{Irreducibility of Primitive Special Cycles in Locally Symmetric Spaces of
Orthogonal Type}
\author{Niklas Ludwig}
\maketitle

\section*{Abstract}
Let $L$ be an even lattice of signature $(n_1,n_2)$ with $n_1\geq k+1,n_2\geq 2$ and let $X$ denote the quotient of the symmetric space associated to $\OG(L\otimes_\ZZ\RR)$ by the discriminant kernel of positive spinor norm. We show that primitive special cycles in $X$ of real codimension at most $kn_2$ are irreducible if $L$ splits off $k+1$ hyperbolic planes locally at every prime.  

\section{Introduction}

Let $L$ be an indefinite even lattice in a rational quadratic space $V$ of signature $(n_1,n_2)$ and let $X$ denote the quotient of the symmetric space associated to $\OG(V\otimes_\QQ \RR)$ by the discriminant kernel in  $\OG(L)\cap\OP(V\otimes_\QQ \RR)$. Kudla and Millson considered \emph{special cycles} in $X$ that are certain linear combinations of subsymmetric spaces associated to $k$-tuples in $V$. Li and Ye show in \cite{ly2026} that the span of the special cycles coincides with the span of adelic special cycles in the rational Chow group in signature $(n,2)$. They ask whether the classes of individual irreducible components belong to the same span. We answer this affirmatively in complex codimension up to $k$ under the hypothesis below.

In \cite{bz2026} so called primitive special cycles are introduced and it is shown that they are irreducible for $k\leq2$ under certain conditions on the lattice $L$. We shall extend this result to $k>2$ under the condition that $L_p\coloneqq L\otimes_\ZZ \ZZ_p$ splits off sufficiently many hyperbolic planes for every prime $p$. More precisely, our main result is the following.

\begin{theorem}
Let $L$ be an even lattice of signature $(n_1,n_2)$ with $n_1\geq k+1,n_2\geq 2$ for an integer $k\geq1$. Suppose $L_p$ splits off $k+1$ hyperbolic planes for every prime $p$.
Then the primitive special cycles in $X$ of codimension at most $kn_2$ are irreducible.
\end{theorem}

We shall also give sufficient criteria for $L$ splitting off hyperbolic planes locally, as follows. 
Let $l(D_L)$ denote the minimal number of generators of the discriminant group $D_L$ of $L$. Most notably, we have the following result.

\begin{theorem}
Suppose $L$ is an even lattice and $l(D_L)+1+2k\leq n_1+n_2$ for an integer $k\geq 1$. Then $L_p$ splits off $k$ hyperbolic planes for every prime $p$.
\end{theorem}

It is well-known that if $L$ is unimodular, then $L_p$ is an orthogonal sum of hyperbolic planes for every prime $p$. For instance, the lattice $U\oplus U\oplus E_8^{\oplus j}$ splits off $2+4j$ hyperbolic planes locally. Here $U$ denotes the hyperbolic plane lattice and $E_8$ is the unique positive definite even unimodular lattic of rank $8$.

\subsection{Acknowledgement}
I would like to thank Jan Bruinier for suggesting the topic of the master's thesis on which this paper is based on, as well as for his guidance and valuable advices. I am also very grateful to Riccardo Zuffetti for many insightful discussions and helpful comments.
The author acknowledges support by Deutsche
Forschungsgemeinschaft (DFG) through the Collaborative Research Centre TRR
326 “Geometry and Arithmetic of Uniformized Structures”, project number
444845124.

\subsection{Locally symmetric spaces}

Let $\OG(V)$ be the rational orthogonal group and let $\OG(L)$ denote the
stabilizer of $L$. Write $V_K=V\otimes_\QQ K$ where $K=\QQ,\QQ_p$ or $\RR$. The spinor norm
$\theta\colon\OG(V_K)\to K^\times/(K^\times)^2$ sends reflection in an
anisotropic vector $v$ to $[q(v)]$, where $[a]=a(K^\times)^2$ denotes
the square class of $a\in K^\times$ \cite[p.134]{b2008}.
Write $\OU(V_K)=\ker\theta$, $\SOU(V_K)=\OU(V_K)\cap\SOG(V_K)$ and let $\OP(V)$ denote the subgroup of positive spinor norm. The symmetric space associated to $\OG(V_\RR)$ can be realized as the Grassmannian
$$ \Dm \coloneqq \{W\subset V_\RR: \dim W = n_2\text{ and } q_{|W}<0\}, $$
which has dimension $n_1n_2$. We denote by $\OT(L)\coloneqq\ker(\OG(L)\to \Aut(D_L))$ the \emph{discriminant kernel}. Let $\Gamma\subseteq\OP(V)\cap\OT(L)$ be a finite-index subgroup and define the orbifold quotient $X\coloneqq \Gamma\backslash\Dm$.
If $\Gamma$ is torsion-free, $X$ is a differentiable manifold.

Let $x=(x_1,\dots x_k)\in V^k$ such that $V_x\coloneqq \spn\{x_1,\dots,x_k\}\subset V$ is positive-definite
and set $r(x)\coloneqq \dim V_x$. The subspaces
$$ Z_x\coloneqq \{W\in \Dm: V_x\perp W\} $$
are sub-symmetric spaces associated to $\OG(V_x^\perp)$
and are thus of dimension $(n_1-r(x))n_2$.

For $T\in\Sym_k(\QQ)$ positive semi-definite, we set 
$$\Omega_k(T)\coloneqq \{x\in V^k: q(x)=T \text{ and } q_{|V_x}>0\}.$$
Kudla and Millson \cite{km1990} constructed \emph{special cycles} of the form 
$$ Z(T,\alpha) \coloneqq \Gamma\backslash \sum_{\substack{x\in\alpha+L^k \\x\in \Omega_k(T)}} Z_x$$
of codimension $(\rk T)n_2$ where $\alpha \in D_L^k$ and $T\in\Sym_k(\QQ)$ positive semi-definite. 
As shown in \Cite{bz2026}, the special cycles can be written as an effective sum of \emph{primitive special cycles} which are given by
$$ \ZP(T,\alpha) \coloneqq \Gamma\backslash \sum_{\substack{x\in\Prm_k(L,\alpha)\\x\in \Omega_k(T)}} Z_x$$
where $\Prm_k(L,\alpha)$ denotes the set of vectors in $\alpha+L^k$ spanning a primitive sublattice of $L'$. They further show using M\"obius inversion, that the primitive special cycles can be retrieved as a linear combination of the special cycles.

We are interested in the irreducibility of the primitive special cycles under certain conditions on $L$. Here we call a (primitive) special cycle in $X$ \emph{irreducible} if its support is the image of a single $Z_x$ in $X$ or the zero cycle.
If $L$ splits off two hyperbolic planes and $k=1$, the Eichler criterion gives a sufficient condition for the irreducibility of the primitive special cycles. A generalization to $k=2$ has been presented in \Cite{bz2026} under the assumption that $L$ splits off $U\oplus U\oplus E_8$.
We give the following sufficient criterion for the irreducibility of the primitive special cycles.

\begin{theorem}
Let $L$ be an even lattice of signature $(n_1,n_2)$ and $n_1\geq k+1,n_2\geq 2$. Suppose that $\SOU(V)\cap\OT(L)\subseteq\Gamma$ and that $L_p$ splits off $k+1$ hyperbolic planes for every prime $p$.
The primitive special cycles in $X$ of codimension at most $kn_2$ are irreducible.
\end{theorem}
\begin{proof}
Let $\ZP(T,\alpha)$ be a primitive special cycle with $\alpha\in D_L^k$ and $T\in\Sym_k(\QQ)$ positive definite. For arbitrary $x,y\in \Prm_k(L,\alpha)\cap \Omega_k(T)$, we set $E_1\coloneqq V_x\cap L'$ and $E_2\coloneqq V_y\cap L'$.
As $q(x)=q(y)$ and $x,y\in \alpha+L^k$, there exists an isomorphism $\phi\colon E_1\to E_2$ with $\phi(x)-x\in L$ for $x\in E_1$.
\Cref{orbit_prim} implies that there exists $\gamma\in\SOU(V)\cap\OT(L)$ with $\gamma(E_1)=E_2$.    
\end{proof}

\section{A Generalized Eichler criterion}

Let $k\in\NN$.
In order to show the irreducibility of the primitive special cycles for $k>2$, a generalization of the Eichler criterion is needed.

We denote the hyperbolic plane lattice over $\ZZ$ and $\ZZ_p$ by $U$ and $U_p$ respectively. 
First, we want to study global conditions that imply that a lattice splits off hyperbolic planes locally. 
As a starting point, there is the following well-known result.

\begin{proposition} \label{un_loc_hyper}
Suppose $L$ is an even unimodular lattice. Then $L_p$ is an orthogonal sum of hyperbolic planes for every prime $p$.
\end{proposition}
\begin{proof}
Let $(n_1,n_2)$ denote the signature of $L$.
As $d(L_2)=(-1)^{n_2}$, \cite[p.78, Theorem 5.2.5]{k1993} implies that $L_2$ is isomorphic to an orthogonal sum of $m$ hyperbolic planes 
where $m\coloneqq \frac{n_1+n_2}{2}$. Since $L$ is unimodular, we have $8\mid n_1-n_2$.
\cite[p.76, Theorem 5.2.4]{k1993} yields that $L_p\cong U_p^{\oplus m}$ for every prime $p\neq 2$. 
\end{proof}

The next result can be seen as a generalization to \cite[Satz 15.8]{k2002}.

\begin{lemma} \label{loc_hyper_split}
Let $L$ be an even $p$-adic lattice and suppose $3+l(D_L)\leq \rk L$. Then $L$ splits off a hyperbolic plane.
\end{lemma}
\begin{proof}
First, we suppose that $p\neq 2$. Then $L$ has an orthogonal basis $x_1,\dots,x_n\in L$ and 
$$ D_L \cong \prod_{i=1}^n \ZZ_p/(x_i,x_i)\ZZ_p. $$
We set $S\coloneqq \{1\leq i\leq n: q(x_i)\in\ZZ_p^\times\}$ and the assumption yields that $|S|\geq 3$.
By \cite[Satz 15.8]{k2002}, the unimodular sublattice $\spn\{x_i: i\in S\}$ splits off a hyperbolic plane.
Secondly, we assume $p=2$ and consider an orthogonal decomposition $L=L_1\oplus\cdots\oplus L_m$ where
the $L_i$ are sublattices of rank $1$ and $2$. Then we claim $L_i$ is not unimodular if and only if $\rk L_i = l(D_{L_i})$.
Indeed, this is valid if $\dim L_i = 1$. If $\dim L_i=2$ and $L_i$ is not unimodular, we choose a basis $x_1,x_2\in L_i$. Then
$$ \det G(x_1,x_2) = \begin{vmatrix} (x_1,x_1) & (x_1,x_2) \\ (x_2,x_1) & (x_2,x_2) \end{vmatrix} \in 2\ZZ_2 $$
and thus $(x_1,x_2)\in 2\ZZ_2$. Therefore $(x,y)\in 2\ZZ_2$ for all $x,y\in L_i$ and thus $l(D_{L_i})=2$.
Let $S\coloneqq \{1\leq i\leq m: L_i=L_i'\}$ and consider the sublattice $E\coloneqq \bigoplus_{i\notin S} L_i$.
Then $\rk E = l(D_E) = l(D_L) \leq \rk L - 3$ and thus 
the unimodular sublattice $F\coloneqq \bigoplus_{i\in S} L_i$ has rank $\rk F\geq 3$. \cite[Satz 15.8]{k2002} yields that $F$ splits off a hyperbolic plane.
\end{proof}

\begin{proposition}
Suppose $L$ is an even lattice and $l(D_L)+1+2k\leq \rk L$. Then $L_p$ splits off $k$ hyperbolic planes for every prime $p$.
\end{proposition}
\begin{proof}
We show by induction over $k$ that $L_p$ splits off $k$ hyperbolic planes.
For $k=1$, the statement follows by \Cref{loc_hyper_split}. Now, assume $k>1$. 
Again \Cref{loc_hyper_split} yields that $L_p=U_p\oplus E$ for some even $p$-adic lattice $E$. The induction hypothesis implies that $E$
splits off $k-1$ hyperbolic planes.
\end{proof}

The main results of this section are based on the principle that splitting off sufficiently many hyperbolic planes locally in a lattice ensures that one can find a hyperbolic plane orthogonal to a given subspace.

\begin{lemma} \label{split_hyper_gl}
Let $L$ be an even lattice of signature $(n_1,n_2)$ with $n_1\geq k, n_2\geq 1$ such that $L_p$ splits off $k$ hyperbolic planes for every prime $p$. If $W\subset L\otimes_\ZZ \QQ$ is a positive definite subspace of dimension at most $k-1$, 
there exists a hyperbolic plane $U\subset L$ with $U\perp W$.  
\end{lemma}
\begin{proof}
By assumption, $L_p$ splits off $U_p^{\oplus k}$. Let $P\colon L_p\otimes_{\ZZ_p} \QQ_p \to U_p^{\oplus k}\otimes_{\ZZ_p} \QQ_p$ denote the orthogonal projection and set $E\coloneqq P(W_p)\cap U_p^{\oplus k}$. To simplify notation, we can write $U_p^{\oplus k}=F\oplus U_p^{\oplus k-1}$
with $F=U_p$. We choose a primitive embedding $\tau\colon E\hookrightarrow U_p^{\oplus k-1}$.
By \cite[Corollary 1.9.7]{n1979}, there exists a $\gamma\in\OG(F\oplus U_p^{\oplus k-1})$ with $\gamma(\im\tau )=E$.
We set $\iota\coloneqq\gamma_{|F}$ and obtain $\iota(F) \perp E$ as $F\perp \im\tau$. 
In particular, $U_p\cong \iota(F) \perp W_p$. As $W^\perp$ is indefinite and $U_p\subseteq L_p\cap W_p^\perp$, 
\cite[Theorem 1.14.2]{n1979} yields that the genus of $L\cap W^\perp$ only contains one class. 
Therefore, there exists an embedding $U\hookrightarrow L\cap W^\perp$ by \cite[p.134, Theorem 6.1.3]{k1993}.
\end{proof}

Suppose $f\in L$ is isotropic and $v\in L$ with $f\perp v$.
Then the isometries $E(f,v):L\to L,$
$$ E(f,v)(x) = x + (x,f)v-(x,v)f-(x,f)q(v)f $$
are called \emph{Eichler transformations}.
We note that $E(f,v)\in\OT(L)\cap\SOU(V)$.

The proof of the following lemma is based on \cite[Lemma 4.14]{bz2026} and uses Eichler transformations.

\begin{lemma} \label{eichler_trans}
Suppose $L=U\oplus E$ for some even lattice $E$. Let $a\in (E')^k$
and set $b\coloneqq (a_1,\dots,a_{k-1},b_k)$ where $b_k\in a_k+E$.
If the entries of each tuple $a$ and $b$ are linearly independent, span a primitive sublattice of $L'$ and $q(a)=q(b)$, 
then there exists $\gamma\in\SOU(V)\cap\OT(L)$ fixing $a_1,\dots,a_{k-1}$ such that $\gamma(a_k)=b_k$.
\end{lemma}
\begin{proof}
As $a_1,\dots,a_k$ span a rank $k$ primitive sublattice of $E'$, they can be extended to a basis of $E'$. By considering the corresponding dual basis,
we find $a'\in E$ satisfying $(a_i,a')=0$ for $i=1,\dots,k-1$ and $(a_k,a')=1$. In the same way, there exists
$b'\in E$ satisfying $(b_i,b')=0$ for $i=1,\dots,k-1$ and $(b_k,b')=1$. Let $f_1,f_2$ be the standard basis of $U$. 
A straightforward calculation shows that
\begin{align*}
&\E(f_1,a')(a_i)=\begin{cases} a_i, &i<k \\ a_k-f_1, &i=k\end{cases} \quad\text{and}\quad \E(f_1,b')(b_i)=\begin{cases} b_i, &i<k \\ b_k-f_1, &i=k.\end{cases}
\end{align*}
Furthermore, we have $\E(f_2,a_k-b_k)(a_k-f_1) = b_k-f_1 $ and $\E(f_2,a_k-b_k)(a_i) = b_i$ for $i=1,\dots,k-1$
by using the assumption that $q(a)=q(b)$.    
\end{proof}

\begin{proposition} \label{orbit_prim}
Suppose $L$ is an even lattice of signature $(n_1,n_2)$ with $n_1\geq k+1,n_2\geq 2$ such that $L_p$ splits off $k+1$ hyperbolic planes for every prime $p$.
If $E_1,E_2\subset L'$ are positive-definite primitive sublattices of rank $k$ and $\phi\colon E_1\to E_2$ an isomorphism with $\phi(x)-x\in L$ for $x\in E_1$, then there exists $\gamma\in\SOU(V)\cap\OT(L)$ such that $\gamma_{|E_1} = \phi$.
\end{proposition}
\begin{proof}
We prove the statement by induction over $k=\rk E_1$. For $k=1$, the statement is \cite[Lemma 4.4]{fr2000}. Now, let $k>1$. We choose a decomposition $F\oplus\ZZ a$ of $E_1$ with $a\in E_1\setminus\{0\}$. By induction, we can assume that 
$\phi_{|F}=\id_F$. 
Let $E\coloneqq F\cap L$. Due to \Cref{split_hyper_gl},
there exists a hyperbolic plane $U\subset L$ with $U\perp E$. Let $L_1\coloneqq U^\perp$. \cite[p.83, Lemma 5.3.1]{k1993} yields that $L_1$ splits off $k$ hyperbolic planes locally.
By applying \Cref{split_hyper_gl} again, we find $\tilde{U}\subset L_1$ with $\tilde{U}\perp E$. Let $P\colon L'\to U\oplus\tilde{U}$ denote the orthogonal projection.
\cite[Lemma 4.4]{fr2000} implies that there exists $\psi\in\SOU(U\oplus \tilde{U})$
such that $\psi(P(a))\perp U$. Therefore, we can further assume that $a\in L_1'$ and similarly $\phi(a)\in L_1'$. Then, we can apply \Cref{eichler_trans}.
\end{proof}

\begin{remark}
Suppose $L$ is an even unimodular lattice. As shown in \cite[Theorem 4.8]{w1964}, the group $\OG(L\oplus U^3)$ is generated by reflections.
Hence, $\theta(\OG(L))\subseteq\{[1],[-1]\}$. In particular, $\OP(L)=\OU(L)$ in this case.
\end{remark}

Suppose the signature of $L$ is $(n,2)$. Splitting off $k+1$ local hyperbolic planes forces $n+2\geq 2(k+1)$. For even $n$
and a lattice containing $U$, we show that the weaker requirement of
extending automorphisms of primitive sublattices already forces $n\geq2k$.
We write $\tau_v$ for reflection in an anisotropic vector $v$, so
$\tau_v(w)=w-(w,v)v/q(v)$ and $\theta(\tau_v)=[q(v)]$.
Thus $[1]$ denotes the trivial square class. For an odd prime $p$,
$\bigl(\frac{a}{p}\bigr)$ is the Legendre symbol: it is $0$ if $p\mid a$,
and otherwise is $1$ or $-1$ according as $a$ is a square or a nonsquare
in $\FF_p^\times$.

\begin{lemma} \label{lift_sublattice}
Let $L$ be a lattice of signature $(n_1,n_2)$ and $E\subseteq L_p$ be a (primitive) sublattice for some prime $p$.
If $\rk E\leq n_1$, then there exists a positive definite (primitive) sublattice $F\subseteq L$ such that $F_p \cong E$.
\end{lemma}
\begin{proof}
By rescaling $L$ if necessary, we can assume that $L$ is even.
Let $k\coloneqq \rk E$ and $x\in E^k$ be a basis of $E$. Since bilinear forms are continuous, there exists $j\in\NN$
such that 
$$G(x)-G(y)\in 4d(E)^2p\M_k(\ZZ_p) \quad\text{if } x-y\in p^jL_p^k.$$
As $L\hookrightarrow L_p$ is dense, we can find $w\in L^k$ such that $x-w\in p^jL_p^k$.
Let $z\in L^k$ be a basis of a positive definite sublattice of $L$ of rank $k$. Then, there exists $i>j$
such that the entries of $y\coloneqq w+p^i z$ span a positive definite sublattice of $L$.
Furthermore, we obtain 
$$x-y=x-w-p^i z\in p^jL_p^k.$$
Consequently, $G(x)-G(y)\in 4d(E)^2p\M_k(\ZZ_p)$. We set 
$$F\coloneqq\spn_\ZZ\{y_1,\dots,y_k\}\subseteq L.$$
Then \cite[Corollary 2.4]{oh2026} yields that $F_p\cong E$.
Now, we assume that $E\subseteq L_p$ is primitive and let $\tilde{F}\coloneqq L\cap W$ where $W\coloneqq\spn_\QQ F$. 
Since $x-y\in pL_p^k$, the lattice $F_p$ is primitive in $L_p$. Therefore, $\tilde{F}_p = L_p\cap W_p = F_p \cong E$.
\end{proof}

\begin{lemma} \label{spinor_lattice}
Suppose $L$ is a lattice. Then $\theta(\OG(L))\subseteq\{[\pm a]: a\mid 2d(L)\}$.
\end{lemma}
\begin{proof}
Let $\phi\in\OG(L)$. Then $\theta(\phi)=\pm\prod_i p_i(\QQ^\times)^2$ for some primes $p_1,\dots p_j$. Further, let $p\neq 2$ be a prime with $p\nmid d(L)$ and $\psi$ be the image 
of $\phi$ under $\OG(L)\to\OG(L_p)$. As $L_p$ is unimodular, we obtain $\theta(\psi)\in\ZZ_p^\times(\QQ_p^\times)^2$.
Therefore, $p\neq p_i$ for $1\leq i\leq j$.
\end{proof}

\begin{lemma} \label{inf_prim_sqr}
Let $k\in\NN$. Then, there exist infinitely many primes $p$ such that $8\mid p+1$ and 
$\genfrac{(}{)}{}{1}{a}{p}=1$ for every $a\in\NN$ with $a\mid k$.
\end{lemma}
\begin{proof}
By Dirichlet's theorem, there exist infinitely many primes $p$ such that $8k\mid p+1$. 
Let us fix such a prime $p$.
Since $8\mid p+1$, the prime $2$ is a square in $\FF_p$.
Let $q\neq2$ be a prime such that $q\mid k$. By the law of quadratic reciprocity, we obtain 
\begin{equation*}
\genfrac{(}{)}{}{1}{q}{p} = \genfrac{(}{)}{}{1}{p}{q} (-1)^{\frac{p-1}{2}\frac{q-1}{2}} = \genfrac{(}{)}{}{1}{-1}{q} (-1)^{\frac{q-1}{2}} = 1. \qedhere   
\end{equation*}
\end{proof}

\begin{lemma} \label{exp_spinor} 
Suppose $p\neq 2$. Then the $p$-adic lattice
$$ E\coloneqq \bigoplus_{i=0}^k \langle 2p^i \rangle $$
satisfies $\theta(\OG(E)) \subseteq \{[1], [p]\}$.
\end{lemma}
\begin{proof}
By \cite[p.81, Theorem 5.3.3]{k1993}, the group $\OG(E)$ is generated by reflections.
Therefore, it is enough to show $\theta(\tau_x)\in \{[1],[p]\}$
for any primitive vector $x\in E$ with $\tau_x\in\OG(E)$.
Let $I$ be the set of those $i\in\{1,\dots,k+1\}$ such that $p^{i-1} x_i^2 \mid (x,x)$.
We show that $x_i\in\ZZ_p^\times$ for $i\in I$. Indeed, suppose that there exists $i\in I$
such that $x_i\in p^m\ZZ_p^\times$ for some $m\in\NN$. Then $p^{i-1+2m}\mid (x,x)$ and a calculation shows
$$e_i - \tau_x(e_i) = \frac{1}{cp^m}x$$
for some $c\in\ZZ_p$. This is a contradiction as $\tau_x\in\OG(E)$ and $x$ is primitive in $E$. Therefore, $x_i\in\ZZ_p^\times$ for $i\in I$.
We let $j\coloneqq \min I$ and see that $q(x)\in p^{j-1}(\ZZ_p^\times)^2$.
\end{proof}

\begin{proposition}
Suppose $L=U\oplus F$ for some even lattice $F$ of signature $(n-1,1)$ with $n\in\NN$ even. If for any positive definite primitive sublattice $E\subseteq L'$ of rank at most $k$ and
$\phi\in\OG(E)$ with $\phi(x)-x\in L$ for $x\in E$, there exists $\gamma\in\OP(V)\cap\OT(L)$ such that $\gamma_{|E}=\phi$,
then $n \geq 2k$.
\end{proposition}
\begin{proof}
We assume by contradiction that $n<2k$. As $n$ is even, this yields $j\coloneqq\frac{n}{2}\leq k-1$. 
By \Cref{inf_prim_sqr}, there exists a prime $p$ such that $4\mid p-3$ and $\genfrac{(}{)}{}{1}{a}{p}=1$ for every $a\mid 2d(L)$.
In particular, $L_p$ is unimodular, $\genfrac{(}{)}{}{1}{-1}{p}=-1$ and \Cref{spinor_lattice} implies that $\OP(L)\subseteq\OU(L_p)$.
The lattice $F_p$ primitively represents the lattice 
$$ \tilde{E} \coloneqq \langle 2(-1)^j \rangle \oplus \bigoplus_{i=2}^j \langle -2p^i \rangle. $$
Thus, \Cref{lift_sublattice} implies that there exists a positive definite primitive sublattice $\tilde{F}\subseteq F'$ such that $ \tilde{F}_p \cong \tilde{E}$.
We fix some primitive embedding $\langle 2(p-1)p \rangle\hookrightarrow U$ which induces a primitive embedding of $E\coloneqq \langle 2(p-1)p \rangle \oplus \tilde{F}$ into 
$L'=U \oplus F'$. Further, we set $\phi\coloneqq \tau_{e_1}\in\OG(E)$ where $e_1$ denotes the standard basis vector of $\langle 2(p-1)p \rangle$.
By assumption, there exists $\gamma\in\OP(L)$ such that $\gamma_{|E}=\phi$. Let $\psi$ be the image of $\gamma$ under $\OP(L)\to\OU(L_p)$.
Since $\theta(\psi_{|E_p})=[-p]$, we get $\theta(\psi_{|E_p^\perp})= [-p]$.
\cite[Corollary 1.6.2]{n1979} and \cite[Corollary 1.9.3]{n1979} yield that
$$ E_p^\perp \cong \bigoplus_{i=0}^j \langle 2p^i \rangle.  $$
Therefore $\theta(\OG(E_p^\perp)) \subseteq \{[1],[p]\}$ by \Cref{exp_spinor}, and we obtain a contradiction.
\end{proof}

The following corollary follows also from \cite[Proposition 1.14.1]{n1979}.

\begin{corollary}
Suppose $L$ is an even unimodular lattice of signature $(n,2)$. If for any positive definite primitive sublattice $E\subseteq L$ of rank at most $k$ and
$\phi\in\OG(E)$, there exists $\gamma\in \OP(L)$ such that $\gamma_{|E}=\phi$,
then $n \geq 2k$.
\end{corollary}

\printbibliography

\end{document}